\magnification=\magstep1
\catcode`\@=11
\input diagrams.sty

\voffset=.75in
\hoffset=.4in

\font\tengoth=eufm10  \font\fivegoth=eufm5
\font\sevengoth=eufm7
\newfam\gothfam  \scriptscriptfont\gothfam=\fivegoth
\textfont\gothfam=\tengoth \scriptfont\gothfam=\sevengoth
\def\goth{\fam\gothfam\tengoth}

\def\dd{{\partial}}

\def\dim{{\rm dim}}

\def\GL{{\rm GL}}

\def\ZZ{{\bf Z}}
\def\FF{{\bf F}}
\def\HH{{\bf H}}
\def\coker{{\rm coker}}

\def\Sym{{\rm Sym}}
\def\Grass{{\rm Grass}}

\def\g{{\goth g}}
\def\gl{{\goth gl}}
\def\d{{\bf d}}
\def\e{{\bf e}}

\def\SS{{\cal S}}

\newcount\cols

{\catcode`,=\active\catcode`|=\active
\gdef\Young(#1){\hbox{$\vcenter
{\mathcode`,="8000\mathcode`|="8000
\def,{\global\advance\cols by 1 &}%
\def|{\cr       
      \multispan{\the\cols}\hrulefill\cr 
      &\global\cols=2 }%
  \offinterlineskip\everycr{}\tabskip=0pt  
  \dimen0=\ht\strutbox \advance\dimen0 by \dp\strutbox    
    \halign
    {\vrule height \ht\strutbox depth \dp\strutbox##
      &&\hbox to \dimen0{\hss$##$\hss}\vrule\cr
     \noalign{\hrule}&\global\cols=2 #1\crcr
     \multispan{\the\cols}\hrulefill\cr%
    }
}$}}
}

\def\sqr#1#2{{\vcenter{\vbox{\hrule height.#2pt
\hbox{\vrule width.#2pt height#1pt \kern#1pt
\vrule width.#2pt}
\hrule height.#2pt}}}}

\def\square{\mathchoice\sqr34\sqr34\sqr{2.1}3\sqr{1.5}3}

\rightline{October 23, 2009}
\bigskip
\centerline {\bf The Existence of Pure Free Resolutions}
\bigskip
\centerline{  David Eisenbud, Gunnar Fl\o ystad and Jerzy Weyman}
\bigskip
\centerline{\it Dedicated to J\"urgen Herzog, on the occasion of his sixty-fifth birthday}
\vskip 1cm

\noindent{\bf Abstract:}
Let $A=K[x_1,\dots, x_m]$ be a polynomial ring in $m$ variables and let
$\d=(d_0<\cdots< d_m)$ be a strictly increasing sequence of $m+1$ integers.
Boij and S\"oderberg conjectured the existence 
graded $A$-modules $M$ of finite length 
having {\it pure free resolution} of type $\d$ in the sense
that for $i=0,\dots,m$ the $i$-th syzygy module of $M$ has generators only in
degree $d_i$.
This paper provides a construction, in characteristic zero, of
modules with this property that are also $GL(m)$-equivariant. Moreover, the construction
works over rings of the form $A\otimes_K B$ where $A$ is a polynomial
ring as above and $B$ is an exterior algebra.

\bigskip
\noindent{\bf Introduction}
\medskip

Let $\d=(d_0<\cdots< d_m)$ be a strictly increasing sequence of integers.
In their remarkable paper  [2006],  Boij and S\"oderberg conjectured the existence 
of a graded module $M$ of finite length over every polynomial
ring $A=K[x_1,\dots, x_m]$ whose minimal free resolution has the form
$$
0\to A^{\beta_m}(-d_m) \to\cdots\to  A^{\beta_0}(-d_0).
$$
Such a resolution is said to be {\it pure\/} of type $\d$.
This paper provides a construction, in characteristic zero, that gives more:
a $GL(m)$-equivariant module with a pure resolution.
In addition we give a another construction of pure resolutions,
of modules supported on determinantal varieties.

The constructions in this paper gave the first proof that the 
Boij-S\"oderberg existence conjecture was correct over any field.
A subsequent paper, by Eisenbud-Schreyer [2007] has verified the conjecture
(along with the other conjectures in the Boij-S\"oderberg paper)
over arbitrary fields, but with (generally) much worse bounds on the ranks of
the modules constructed, and with no $GL(m)$-equivariance. 
 Also, Sam and Weyman [2009] have provided
a more direct proof that the modules we construct actually have pure resolutions.

One of the 
major questions left open by both these papers
is the nature of the semigroup of possible degrees of modules having pure
resolutions with a given degree sequence. There is a unique minimum
possibility, determined solely by integrality considerations, and it is easy to see that only
integral multiples of this minimum can occur. It is known from many examples,
that not all occur. We make some conjectures in this direction,
supported by the examples in this paper, in section 6. In particular, 
Conjecture 6.1 asserts  that every sufficiently high multiple that is integral
{\it does\/} occur. The paper of Eisenbud-Schreyer [2007] is 
complementary to this one, in that putting together
the two equivariant examples produced here with the example produced there
shows that this conjecture is true for many degree sequences. 

We are able to make our constructions in a more general context
than the polynomial rings that occur in Boij-S\"oderberg [2006] or
Eisenbud-Schreyer [2007]: we work over a free
strictly commutative $\ZZ/2$-graded algebra, which specialize to 
the necessary pure resolutions in the polynomial algebra case, and also give 
pure resolutions over
exterior algebras. This suggest that there may be a stronger version
of the Boij-S\"oderberg conjectures
 Eisenbud-Schreyer theory addressing resolutions over
 the exterior algebra as well.
 
Our constructions make use of Schur functors, $\ZZ/2$-graded Schur functors and Bott's 
Theorem. This last is what limits
our method to characteristic 0.
It remains an open
question whether such examples exist in characteristic $p>0$.


Let $V=V_0\oplus V_1$ be a $\ZZ/2$-graded vector space with dimension vector $(m, n)$. We recall (see Section 1)  that there exist $\ZZ/2$-graded versions of Schur modules
$$
\SS_\lambda (V)=\bigoplus_{\mu\subset \lambda}
 S_\mu V_0\otimes S_{\lambda'/\mu'}V_1 .
$$
Here $\lambda'$ denotes the conjugate partition to $\lambda$. For example, the
conjugate partition to $(2)$ is $(1,1)$.

We work over  the free $\ZZ/2$-graded algebra 
$$R= Sym (V)=\bigoplus_{i\ge 0} \SS_i (V)= Sym(V_0 )\otimes  \bigwedge{}^\bullet (V_1).$$

Let $\lambda = (\lambda_1 ,\ldots ,\lambda_s ,\dots )$ be a partition and $e_1$ an integer.
We use the convention that $\lambda$ has infinitely many parts, but only finitely many are non-zero. 
Define the integers $e_i$ for $2\le i$, by setting
$$e_i =\lambda_i -\lambda_{i-1}+1,
$$
so that $e_i=1$ for $i\gg 0$. For convenience we set $d_i=\sum_{j=1}^i e_j$, and 
$d_0=0$. We also set $\d=(d_0,d_1,\dots)$. 

Next, define a sequence of partitions $\alpha(\d, i)$ for $i\ge 0$. We set
$\alpha(\d, 0)=\lambda$ and, for $i\geq 1$,
$$
\alpha(\d, i)= (\lambda_1+e_1 ,\ldots ,\lambda_i+e_i ,\lambda_{i+1},\ldots ,\lambda_s , \dots).
$$
The $\ZZ/2$-graded endomorphisms of $V_0 \oplus V_1$ form a $\ZZ/2$-graded Lie algebra
$\gl(V)$. We define a  complex of free $\gl(V)$-equivariant $R$-modules, with the terms
$$
\eqalign{
\FF (\d )_0&=\SS_\lambda V\otimes R, \cr
\FF (\d )_i &=\SS_{\alpha(\d, i)} V\otimes R(-e_1-\ldots-e_i )\ {\rm for}\ i\ge 1.
}
$$
Though this complex is in general infinite, it becomes finite of length at most
$\dim V_0$ in the case when $V_1=0$.
The differential
$$
\dd_i :\FF(\d )_i\rightarrow \FF(\d)_{i-1}$$
is given on the generators by the $\ZZ/2$-graded Pieri maps (see Section 1)
$$ \SS_{\alpha(\d, i)} V\rightarrow \SS_{\alpha(\d, i-1)}V \otimes \SS_{e_i} V=\SS_{\alpha(\d, i-1)}V\otimes R_i .
$$
We now state the main results of this paper.

\proclaim Theorem 0.1. The complex $\FF(\d )_\bullet$ is an acyclic complex of $\gl(V)$-equivariant,  free $R$-modules. It is also pure, with the $i$-th differential of degree $e_i$.

\medskip
\noindent{\bf Remark.} When $n>0$ the complexes $\FF(\d )_\bullet$ are infinite, but
eventually linear. In fact, if the partition $\lambda$ has $s$ non-zero parts, 
the complex $\FF (\d )_\bullet$ becomes linear after $s$ steps.
\medskip

For the second construction we fix two $\ZZ/2$-graded spaces $V$ and $U$ (with $U$ of dimension vector $(d,e)$) and we work over the symmetric algebra
$S=\Sym (V\otimes U)$. Fix a sequence $\d$ as above. 

We define an infinite complex $\HH(\d)$ of free $\gl (V)\times \gl(U)$-equivariant $S$-modules, with the terms
$$\eqalign{
\HH (\d )_0&=\SS_\lambda V\otimes S, \cr
\HH (\d )_i &= \SS_{\alpha(\d, i)} V\otimes \SS_{d_i}U\otimes S(-e_1-\ldots-e_i )\ {\rm for}\ i\ge 1.
}
$$

The differential
$$
\dd_i :\HH(\d )_i\rightarrow \HH(\d)_{i-1}
$$
is given on the generators by the $\ZZ/2$-graded Pieri maps (see Section 1)
$$
\matrix{
\SS_{\alpha(\d, i)} (V)\otimes  \SS_{d_i}U\cr
\downarrow\cr
 \SS_{\alpha(\d, i-1)}(V)\otimes \SS_{d_{i-1}}U\otimes \SS_{e_i}V\otimes\SS_{e_i} U\cr
 \bigcap\cr
  \SS_{\alpha(\d, i-1)}V\otimes \SS_{d_{i-1}}U\otimes S_i .
}
$$

\proclaim Theorem 0.2. The complex $\HH(\d )_\bullet$ is an acyclic complex of  free $S$-modules that is $\gl (V)\times \gl(U)$-equivariant. It is pure, and the $i$-th differential has degree $e_i$.

The paper is organized as follows. Section 1 is devoted to recalling needed notions from representation theory of $GL(m)$. In \S 2 we briefly review the material needed from representation theory of $\ZZ/2$-graded Lie algebra $\gl(V)$. In \S 3 and \S 4 we prove 
Theorems 0.1 and 0.2 in special cases. In \S 5 we deduce the general cases from
the cases already treated. Further conjectures and open problems are discussed in \S 6.

We are grateful to F.-O. Schreyer for pointing out to us that
certain known complexes constructed
by multilinear algebra (see \S 2)
give examples including all the pure resolutions whose Betti tables have
just two rows with nonzero terms, and thus setting us on the idea
of using Schur functors to construct pure resolutions.

\bigskip
\noindent{\bf \S 1. Cohomology of Homogeneous Bundles on Projective Spaces}
\medskip

For the convenience of the reader we review a few necessary results from representation theory.

We work over a field $K$ of characteristic zero. We denote by $E$ a vector space of dimension $m$ over $K$ (or sometimes a vector bundle
of rank $m$ on an algebraic variety), and write $A = \Sym (E)$.
Here and in the sequel we use the language of
vector bundles but always work with the associated locally free sheaves.

There is a one-to-one correspondence between the irreducible
polynomial representations
of the group $GL(E)$ and the set of {\it partitions\/} 
$\lambda =(\lambda_1 ,\ldots ,\lambda_m )$ with at most $m$ parts. The representation corresponding
 to $\lambda$ will be written $S_\lambda E$, and $\lambda$ is called the
 {\it highest weight\/} for $S_\lambda E$. The construction of these
 representations is functorial; in characteristic 0, for
 example, one may view the representation
 $S_\lambda E$ as the image of $E\otimes\cdots\otimes E$, the tensor
 product with $t=\sum_i\lambda_i$ factors, under the projection map
 defined by a {\it Young symmetrizer\/}, which is a certain element of the 
 group algebra of the symmetric group on $t$ letters that acts by
 permuting the factors of the tensor product (see Fulton and Harris [1991]
 \S 4.1 and \S 6.1.) For this reason the 
 construction extends to the case where $E$ is a vector bundle
 on an arbitrary space, and the proof below will imply that
 the complexes of vector bundles we construct are  resolutions because  acyclicity can be
 proved fiberwise.
  
For example the $d$-th symmetric power of $E$ is
$S_{(d,0,\dots,0)}E$, which we will often denote by $S_d E$.  
Of course $A=\oplus_{d\ge 0} S_d E$. The one-dimensional
representation $\wedge^mE$ corresponds to the weight $(1^m):=(1,1,\dots,1)$. 
For any $\lambda$ and integer $p$ we have 
$$
S_\lambda E\otimes (\wedge^m E)^p=S_{\lambda+(p^m)}
$$
where $\lambda+(p^m) = (\lambda_1+p, \dots, \lambda_m+p)$. Thus we may assume that all the
$\lambda_i$ are non-negative and that $\lambda_m=0$ when this is convenient.
It is useful to visualize $\lambda$ as  a
{\it Young frame\/}, a diagram of boxes in which
the $i$-th column of boxes extends down $\lambda_i$ boxes from a 
given baseline; for example, the Young frame for $\lambda=(4,2,2,1,1)$ is
$$
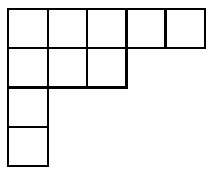
\Young(,,,,|,,||)
$$

There is a general formula giving the decomposition---in 
characteristic 0---of the tensor
product of two representations, called the Littlewood-Richardson Rule.
Here we will only use the simple special case called the {\it Pieri Formula\/}, which gives the decomposition of $S_\lambda E\otimes S_dE$ for any $\lambda$ and $d$.
To express it, we define $|\lambda |:=\sum_{i=1}^m\lambda_i$, and 
we write $\mu\supset\lambda$ if $\mu_i\ge\lambda_i$ for $i=1,\ldots n$. We 
will say that $\mu /\lambda$ is {\it a vertical strip} if $\mu\supset\lambda$ and if
$\mu_i\le\lambda_{i-1}$ for $i=2,\ldots ,n$. In the case where all the 
$\lambda_i$ are non-negative, then $\mu\supset \lambda$
means that the Young frame for $\lambda$ fits into the upper left hand
corner of the Young frame for $\mu$, and that no two boxes of $\mu$
that are outside $\lambda$ lie in the same row.

Since the decomposition is once again given by applying Young symmetrizers,
it works for vector bundles as well.

 \proclaim Theorem 1.1 (Pieri's Formula). If $E$ is
 a vector bundle defined on an algebraic variety of characteristic 0 then
$$
S_\lambda E\otimes S_i E=\oplus_\mu S_\mu E
$$
where the sum is taken over all partitions $\mu\supset \lambda$ such that $|\mu |-|\lambda |=i$ and $\mu/\lambda$ is a vertical strip.

\noindent{\sl Proof.} See Weyman [2003] (2.3.5) or
Fulton-Harris [1991], Appendix A, (A.7). See also MacDonald
[1995], Chapter 1. \hfill$\square$\medskip

The other result from representation theory that we need is a special case of
Borel-Bott-Weyl theory. Let $\Grass (1, E)$ denote the Grassmannian of
1-dimensional subspaces of $E$, which may also be viewed as
a the projective space, 
$$
\Grass(1,E)= {\bf P}(E^*)\cong {\bf P}^{m-1}.
$$
Let  ${\cal R}$ denote the the tautological rank one sub-bundle on $\Grass(1,E)$,
and let ${\cal Q}$ the the quotient bundle, with {\it tautological exact sequence\/}
$$
0\rightarrow {\cal R}\rightarrow E\otimes {\cal O}_{{\bf P}^{m-1}}\rightarrow {\cal Q}\rightarrow 0.
$$

For any sheaf ${\cal G}$ on $\Grass(1,F)$, let $H^i ({\cal G})$ denote the cohomology $H^i ({\bf P}^{m-1}, {\cal G})$. The 
result we need describes this cohomology in the case of an equivariant sheaf
${\cal G}=S_\alpha {\cal Q}\otimes S_u{\cal R}$. To express it we need
two other pieces of notation. For any permutation $\sigma$
we write $l(\sigma)$ for the {\it length\/} of $\sigma$, that is, the minimal number
of transpositions necessary to express $\sigma$ as a product of transpositions.
We write $\rho$ for the partition $\rho=(m-1,m-2,\dots, 1, 0)$.

\proclaim Theorem 1.2 (Bott's Theorem in a special case). 
With notation as above, $H^i(S_\alpha {\cal Q}\otimes S_u{\cal R})$ is nonzero
for at most one index $i$. More precisely, consider the sequence of integers
$(\alpha ,u)+\rho =(\alpha_1+m-1,\ldots ,\alpha_{m-1}+1,u)$.
\item{1)} If the sequence $(\alpha ,u)+\rho$ has a repetition
then the sheaf $S_\alpha {\cal Q}\otimes S_u{\cal R}$ has all cohomology equal to zero.
\item{2)} If the sequence $(\alpha ,u)+\rho$ has no repetitions
then there exists a unique permutation $\sigma$ such that 
$\beta:=\sigma ((\alpha ,u)+\rho)-\rho$
is non-increasing. In this case $S_\alpha {\cal Q}\otimes S_u{\cal R}$
has  only one nonvanishing cohomology group, which is
$$
H^{l(\sigma )}(S_\alpha {\cal Q}\otimes S_u{\cal R} )=S_\beta E.
$$

\noindent{\sl Proof.} The dual form of this result is Weyman [2003], (4.1.9).
The version used here follows by the duality result given in 
Exercise 2.18b in Weyman [2003]. A very short argument for Bott's theorem
may be found in Demazure [1976].
\hfill $\square$

\proclaim Corollary 1.3. Let $\lambda_1\geq\cdots\geq \lambda_{m-1}$ be
a sequence of non-negative integers, and let ${\cal B} = \Sym({\cal Q})$.
\item{a)}If $\lambda_{m-1}=0$, then there is an equivariant
isomorphism of graded $A:=Sym E$-modules
$$ 
H^0 (S_{(\lambda_1 ,\ldots ,\lambda_{m-1})}{\cal Q}\otimes {\cal B})\cong
S_{(\lambda_1 ,\ldots ,\lambda_{m-1},0)}E\otimes A.$$
\item{b)} If $\lambda_{m-1}>0$ then
$H^0 (S_{(\lambda_1 ,\ldots ,\lambda_{m-1})}{\cal Q}\otimes {\cal B})$
has an equivariant
minimal resolution by  free graded $A$-modules of the form
$$
0\rTo 
S_{(\lambda_1 ,\ldots ,\lambda_{m-1},1)}E\otimes A(-1)
\rTo
S_{(\lambda_1 ,\ldots ,\lambda_{m-1},0)}E\otimes A
$$

\noindent{\sl Proof.} From the tautological exact sequence above we derive a resolution
of each $\Sym_d({\cal Q})$, and thus of the graded algebra $\cal B$, which takes the form
$$
0\rightarrow
A(-1)\otimes {\cal R} \rightarrow
 A\otimes {\cal O}_{{\bf P}^{m-1}}\rightarrow 
 {\cal B}\rightarrow 
 0.
$$
We tensor this resolution with $S_{(\lambda_1,\dots,\lambda_{m-1})}{\cal Q}$
and form the long exact sequence in cohomology,
$$\eqalign{
0\rightarrow
A(-1)\otimes H^0 (S_{(\lambda_1 ,\ldots ,\lambda_{m-1})}{\cal Q}\otimes {\cal R})
\rightarrow &
A\otimes H^0 (S_{(\lambda_1 ,\ldots ,\lambda_{m-1})}{\cal Q})
\rightarrow\cr
H^0 (S_{(\lambda_1 ,\ldots ,\lambda_{m-1})}{\cal Q}\otimes {\cal B})
&\rightarrow 
A(-1)\otimes H^1 (S_{(\lambda_1 ,\ldots ,\lambda_{m-1})}{\cal Q}\otimes {\cal R}).
}
$$
 If $\lambda_{m-1} = 0$, Bott's Theorem shows
that all the cohomology of
${\cal R}\otimes S_{(\lambda_1 ,\ldots ,\lambda_{m-1})}{\cal Q}$
vanishes. By Bott's Theorem, $H^0(S_{\lambda} \cal Q)=S_\lambda E$, so we get Part a). If, on the other hand, 
$\lambda_{m-1} > 0$ then Bott's Theorem shows that
the $H^1$ term is zero, and the resulting equivariant short exact sequence
is the one given in Part b). \hfill$\square$

We remark that the use of the complex in b), which is
the push-down of the complex 
$S_\lambda{\cal Q}\otimes\bigwedge^\bullet ({\cal R})$,
is a simple example of the geometric technique described in Weyman [2003].

\bigskip
\noindent{\bf \S 2. $\ZZ/2$-Graded Representation Theory}
\medskip

For the proof of Theorem 0.1 we will use the results
of Berele and Regev [1987] giving the structure of $R$ as
a module over a $\ZZ/2$-graded Lie algebra $\g :=\gl(V)$. For the convenience of the reader
we give a brief sketch of what is needed.
Let $V=V_0\oplus V_1$ be a $\ZZ/2$-graded vector space of dimension $(m,n)$.

The
$\ZZ/2$-graded Lie algebra
$\gl(V)$ is the vector space of $\ZZ/2$-graded
 endomorphisms of $V=V_0\oplus V_1$. Thus
$$\gl(V)= \gl(V)_0\oplus \gl(V)_1,$$
where $\gl(V)_0$ is the set of endomorphisms preserving the grading of $V$
and $\gl(V)_1$ is the set of
endomorphisms of $V$ shifting the grading by 1. Additively
$$\gl(V)_0 = End_K (V_0 )\oplus End_K (V_1 ),$$
$$ \gl(V)_1 =Hom_K (V_0 ,V_1 )\oplus Hom_K (V_1 ,V_0 )$$
The commutator of the pair of homogeneous elements $x,y\in \gl(V)$ is
defined
by the formula
$$[x,y]= xy-(-1)^{deg(x)deg(y)}yx.$$

By
a $\gl(V)$-module we mean a $\ZZ/2$-graded vector space $M=M_0\oplus M_1$
with a bilinear map of $\ZZ/2$-graded vector spaces
$\circ :\gl(V)\times M\rightarrow M$   satisfying the identity
$$
[x,y]\circ m = x\circ (y\circ m)-(-1)^{deg(x)deg(y)} y\circ (x\circ m))
$$
for homogeneous elements $x,y\in \gl(V), m\in M$.

In contrast to the classical theory, not every representation of the
$\ZZ/2$-graded Lie algebra $\gl(V)$ is
semisimple. For example its natural action on mixed tensors $V^{\otimes
k}\otimes V^*{}^{\otimes l}$ is in
general not completely reducible. However, its action on
$V^{\otimes t}$ decomposes just as in the ungraded case:

\proclaim Theorem 2.1 (Berele-Regev 1987). The action of $\gl(V)$ on $V^{\otimes t}$ is
completely reducible for each $t$.
More precisely, the analogue of Schur's double centralizer theorem holds
and the irreducible $\gl(V)$-modules
occurring in the decomposition of $V^{\otimes t}$  are in 1-1
correspondance with irreducible representations of
the symmetric group $\Sigma_t$ on $t$ letters. These irreducibles are the ($\ZZ/2$-graded) Schur
functors
$$\SS_\lambda (V) = e(\lambda )V^{\otimes t}$$
where $e(\lambda )$ is a Young idempotent corresponding to a partition
$\lambda$ in the group ring of the
symmetric group $\Sigma_t$.

(This result is also proven in Section 1 of Eisenbud-Weyman [2003].)
The notation is consistent with the notation above in the sense that the $d$-th
homogeneous component of the ring $\SS (V)$ is $\SS_d (V)$ where $d$ represents the
partition $(d)$ with one part.

Here we use the symbol $\SS_\lambda$ to denote the $\ZZ/2$-graded
version of the Schur functor $S_\lambda$; the latter acts on ungraded vector
spaces. 
The partition $(d)$ with only one part will be denoted simply $d$, so
for example
$\SS_2(V)=S_2(V_0)\oplus (V_0\otimes
V_1)\oplus
\wedge^2V_1$ and similarly
$\bigwedge^{2}V=\SS_{(1,1)}V=\wedge^2V_0\oplus V_0\otimes V_1\oplus S_2(V_1)$.
In each case the decomposition is as representations of the subalgebra
$\gl(V_0)\times \gl(V_1)\subset \gl(V)$. Similar decompositions hold for
all $\SS_d V$ and $\bigwedge^{d}V$. 
The Pieri formula (and the Littlewood-Richardson rule) generalize verbatim to $\ZZ/2$-graded Schur functors:

 \proclaim Proposition 2.2 ($\ZZ/2$-graded Pieri Formula). If $V$ is
 a $\ZZ/2$-graded vector space, and $\lambda$ a partition, we have an isomorphism of $\gl (V)$-modules
$$
\SS_\lambda V\otimes \SS_i V=\bigoplus_\mu \SS_\mu V
$$
where the sum is taken over all partitions $\mu$ such that $|\mu |-|\lambda |=i$ and $\mu/\lambda$ is a vertical strip.

This follows from the results of Berele-Regeve [1987]. 

\bigskip
\noindent{\bf \S 3. First Construction of Pure Resolutions in the Even Case}
\medskip

Let $E$ be an $m$-dimensional vector space, or more generally a rank
$m$ vector bundle on an algebraic variety, over a field of characteristic 0.
Fix a strictly increasing sequence of integers $\d=(d_0 ,d_1 ,\ldots ,d_m)$. 
We will produce a pure acyclic equivariant complex of length $m$ with terms in degrees $d_0 ,\ldots, d_m$. To simplify notation we set
$$
e_0 :=d_0, \qquad e_i: =d_i -d_{i-1}, {\rm for}\ i=1,\ldots ,m,
$$
and we sometimes write $\e=(e_0,\dots,e_m)$ for the sequence
corresponding to $\d$.

We will construct a complex
$$
{\bf F}(\d)_\bullet= {\bf F}(\d)(E)_\bullet :\quad 0\rightarrow F (\d)_m
\rightarrow F(\d )_{m-1}\rightarrow\ldots\rightarrow F(\d )_1\rightarrow F (\d )_0
$$
where $F(\d )_i $ is a free $A=\Sym(E)$-module generated in degree 
$d_i=e_0 +\ldots +e_i$. Let $\lambda = (\lambda_1,\dots,\lambda_m)$ be the partition with parts
$\lambda_i = e_0+ \sum_{i+1}^m (e_j-1)$ (so in particular
$\lambda_m=e_0$). We define a sequence of partitions
$\alpha (\d ,i )$ for $0\le i\le m$ by
$$\eqalign{
\alpha (\d ,0) &= \lambda ,\cr
\alpha (\d, i)&= (\lambda_1 +e_1 ,\lambda_2 +e_2 ,\ldots ,\lambda_i +e_i ,\lambda_{i+1},\ldots ,\lambda_{m}),
}
$$
and set
$$
F (\d)_i:= S_{\alpha (\d,i)}E\otimes A(-e_0 - e_1 -\ldots -e_i).
$$

We could of course have reduced to the case
$d_0=0$; as defined below, the resolutions cases with
$d_0\neq 0$ are obtained from the ones with $d_0=0$
simply by tensoring with the 1-dimensional representation
$(\wedge^m E)^{d_0}$. We will sometimes make the assumption $d_0=0$ for 
simplicity, but we will need the case with $d_0>0$ for 
induction.  

To make it easier to think about these complexes, we give a pictorial 
representation. The following example contains all features of general case.

\medskip
\noindent{\bf Example 3.1.} Take $m=4$ and $\d= (0,4,6,9,11)$,
so that $\e= (0,4,2,3,2)$. Then the partition $\alpha (\d,i)$ (for $0\le i\le 4$) is the subdiagram filled with numbers $\le i$ in the Young diagram
$$
\Young (0,0,0,4|0,0,3,4|0,0,3|0,2,3|1,2|1|1|1).
$$
Here $\alpha(\d ,0)$ is the partition in which difference between the $i$-th and $i+1$'st column equals $e_{i+1}-1$.
We get $\alpha (\d ,i)$ from $\alpha (\d ,i-1)$ by adding $e_i$ boxes to the $i$-th column. 
Notice that for each $i\ge 1$ there is exactly one row
in the diagram above containing
boxes numbered $i$ and $i+1$---
these are the highest box with the number 
$i$ and the lowest box with the number $i+1$.
This is a general phenomenon that makes it possible
for us to define a differential, there is really no
choice about how to construct it because of the following observation:

\medskip

Because $S_{\alpha (\d ,i )} E$ is obtained from
$S_{\alpha (\d ,i-1 )} E$ by adding $e_i$
boxes in one column,
the Pieri Formula implies that it
occurs exactly once in the decomposition of
in $S_{\alpha (\d ,i-1 )} E\otimes S_{e_i}E$
into irreducible $\GL(E)$-modules. Thus there
is a unique (up to scalar) nonzero equivariant map
of $A$-modules
$$
\phi(\d ,i): F(\d)_i \rightarrow F(\d )_{i-1},
$$
and it has degree 0 in the grading coming from $A$
since the generators of 
$F(\d)_i$ have degree $e_i$ more than those of $ F(\d )_{i-1}$.

For any $i\leq m-2$, there is a row of 
$\alpha(\d, i+2)$ containing two more boxes than are present in
$\alpha(\d, i)$. The Pieri Formula thus implies 
that $S_{\alpha (\d ,i )} E$ does not occur in 
$S_{\alpha (\d ,i-2 )} E\otimes S_{e_i +e_{i-1}}E$, so 
$\phi(\d ,i-1)\phi(\d ,i)=0$, so the maps
$\phi(\d ,i)$ make
${\bf F}(\d)_\bullet$ into a complex. This argument
actually shows that any equivariant maps of $\GL(E)$
modules ${\bf F}(\d)_i\to {\bf F}(\d)_{i-1}$ would make 
${\bf F}(\d)_\bullet$ into a complex; and by the construction above,
any nonzero equivariant maps make it into a complex isomorphic to
${\bf F}(\d)_\bullet$. We will use this
uniqueness in the proof below.

Here is our main result in the even case:

\proclaim Theorem 3.2. If $E$ is a vector space of
dimension $m$ over a field of characteristic 0, and 
$\d=(d_0,\dots,d_m)$ is a strictly increasing sequence of 
integers, then
\item{1)} The complex
 $$
{\bf F}(\d)(E)_\bullet:\quad 
0\rTo F{(\d)}_m\rTo^{\phi(\d ,m)}\cdots \rTo F(\d )_{1}\rTo^{\phi(\d ,1)} F(\d)_0\rTo 0
$$ 
is a minimal
graded free resolution, and the generators of $F(\d)_i$ have degree $d_i$.
\item{2)} The module
$
M(\d):= \coker \, \phi(\d ,1).
$
resolved by ${\bf F}(\d)(E)_\bullet$
is equivariant for $\GL(E)$.
As a representation, $M(\d)$ is isomorphic to  the direct sum  of all the 
 irreducible summands of $S_{\alpha (\d,0)}E\otimes \Sym (E)$ 
corresponding to the partitions that do not contain 
$\alpha (\d,1)$. In particular
$M(\d)$ is finite dimensional as a vector space, and is
zero in degrees $\geq \alpha (\d,1)_1$.

\noindent{\bf Remark:\/} If we simply think of each ${\bf F}(\d)_i$ as
a sum of representations, and define $M(\d)$ as the sum of the
representations in part 2), then in the augmented complex consisting
of ${\bf F}(\d)_\bullet$ and $M(\d)$, each irreducible representation that
occurs in one term occurs either in the term before or the term after,
but not both. Moreover, in a given ${\bf F}(\d)_i$ no representation 
occurs more than once. Thus we see that it is 
combinatorially possible that ${\bf F}(\d)_\bullet$ is a resolution of $M(\d)$.
To make this into a proof of Theorem 1, one could first apply the 
Acyclicity Lemma of Peskine and Szpiro [1974], which implies that it is enough
to prove the acyclicity of ${\bf F}(\d)_\bullet$ after replacing
the variables $(x_1,\dots,x_m)$ in $\Sym(E) = K[x_1,\dots, x_m]$ by
$(1,0,\dots,0)$. To finish the proof, one would need to show that
the highest weight vector of each $\GL(x_2, \ldots,x_m)$-representation 
contained in both
$S_{\alpha(\d,i)}$ and $S_{\alpha(\d,i-1)}$ is mapped from the first module
into a nonzero vector in the second. The proof below shows that this
must in fact be true! But we do not at present know how to 
prove it directly.

\bigskip
\noindent{\sl Proof of Theorem 3.2.}
We use induction on $m$ and (in the last part of the proof)
on $d_m-d_0 = \sum_{i\geq 1} e_i.$  If $m=1$ then
the complex has the form
$$
{\bf F}(e_0 ,e_1 )_\bullet :A(-e_1-e_0)\rightarrow A(-e_0)
$$
with the map being the multiplication by $x_1^{e_1}$, and the assertions
are trivial. On the other hand, if
$d_m-d_0=m$, the smallest  possible value, then all the $e_i$ are
1 and the complex ${\bf F}(\d)_\bullet$ is simply the Koszul
complex on the variables in the polynomial ring $A$, so the 
theorem is true in this case as well.

We next show that part 1) of the theorem, for a given $m$,
implies part 2) for that $m$. 
We use Pieri's formula to understand
the $F(\d)_i$, and assume that ${\bf F}(\d)_\bullet$ is a resolution
of $M(\d)$.
Since no $S_\beta E$
occurs with multiplicity more than 1 in a term of the complex,
a representation is present (with multiplicity 1) in
$M(\d)$ if it is present in 
 $F(\d)_0$ but not
 $F(\d)_1$; and it is absent from 
 $M(\d)$ if it is either absent from 
 $F(\d)_0$  or present in
 $F(\d)_0$ and also in 
 $F(\d)_1$ but not in
 $F(\d)_2$.
 
First, if $\beta \not \supset  \alpha(\d, 1)$ then $S_\beta$
cannot occur in $F(\d)_1 = 
A\otimes S_{(\d, 1)} E$, so if $S_\beta$ is
present in $F(\d)_0=A\otimes S_\lambda E$ then it 
is present in $M(\d)$.

Next suppose that $\beta\supset  \alpha(\d, 1)$ and
$S_\beta$ occurs in $A\otimes S_\lambda E$. It is clear from
the Pieri formula that
$S_\beta$ also occurs in $A\otimes S_{\alpha(\d, 1)} E$.
But since $S_\beta$ occurs in $A\otimes S_\lambda E$, and 
$\beta_1\geq \alpha(\d, 1)>\lambda_1$, we must have
$\beta_2\leq \lambda_1$.  It follows that 
$\beta\not\supset \alpha(\d, 2)$. Thus $S_\beta$ does
not occur in $F(\d)_2$, so it is in the image
of $F(\d)_1\to F(\d)_0$,
and thus cannot occur in $M(\d)$, completing the
proof of part 2) based on part 1).
\smallskip

For the inductive step in the proof of part 1),
we consider the sheaf of algebras  ${\cal B}=\Sym ({\cal Q})$ on $\Grass (1, E)\cong 
{\bf P}^{m-1}$, and let
${\bf F}(\d)({\cal Q})$ be the corresponding complex of vector
bundles on $\Grass(1,E)$.
The bundle $\cal Q$ has rank $m-1$, so applying our induction
on the dimension of $E$ to the fibers of the bundle $\cal Q$ at each point,
we see
 that the complex of vector bundles
${\bf F}(e_0 ,\ldots ,e_{m-1})({\cal Q})_\bullet$, and with it the complex
$${\cal F}_{\bullet}:=
{\bf F}(e_0 ,\ldots ,e_{m-1})({\cal Q})_\bullet\otimes (\bigwedge^{m-1}{\cal Q})^{\otimes {e_m -1}}
$$
is acyclic. Its terms are
the Schur functors on $\cal Q$ with highest weights 
$$
\alpha'(\d ,i):=
\alpha(\d ,i)_1,\dots, \alpha(\d ,i)_{m-1}
$$ 
for $0\le i\le m-1$---the same as $\alpha'(\d ,i)$
but with the last part  $\alpha (\d ,i)_m=0$ omitted to
make a partition of length $m-1$. By induction, 
$\cal F_\bullet$ is a resolution of a $\cal B$-module
that we may call $M_{\cal Q}(\d)$, which is a 
 direct sum of finitely many representations, each a Schur functor
 of ${\cal Q}$.

Next consider the complex obtained from $\cal F$ by taking global sections,
$$
H^0({\cal F}_{\bullet}): 0\rightarrow H^0({\cal F}_{m-1})
\rightarrow\cdots\rightarrow
H^0({\cal F}_{0}).
$$
By Bott's Theorem,  $H^j({\cal F}_i)=0$ for all $i$ and all $j>0$.
Breaking the complex $\cal F$ into short exact sequences, one sees
by induction that this implies the acyclicity of the complex
$H^0({\cal F}_{\bullet})$, and this is a resolution of the $A$-module
$H^0(M_{\cal Q}(\d))$.

By the Corollary 1.3,
 each term 
 $$
 H^0({\cal F}_i) = H^0 \biggl({\bf F}(e_0 ,\ldots ,e_{m-1})({\cal Q})_i \otimes (\bigwedge^{m-1}{\cal Q})^{\otimes {e_m -1}})\biggr)
 $$ 
 of $H^0({\cal F}_{\bullet})$ is either free or has a free resolution of length $1$. 
 We distinguish these two cases. The reader will find an explicit example
 for each of these cases in Example 3.3 and Example 3.4 below, and
 it may be helpful to consider the pictures there while reading the following.
 
{\bf Case 1)} $e_{m}=1$. In this case each 
${\bf F}(e_0 ,\ldots ,e_{m-1})({\cal Q})_i \otimes (\bigwedge^{m-1}{\cal Q})^{\otimes {e_m -1}}$
for $i\leq m-2$ satisfies the conditions of 
Part a) of the Corollary to Bott's Theorem. Thus 
the modules $H^0({\cal F}_i)$ for $i<m-1$ are free, and are the same as
those of ${\bf F}(\d)_\bullet$. By part b), on the other hand
the last term $H^0({\cal F}_{m-1})$
has homological dimension 1, and we see that the terms of its resolution furnish the remaining two terms of ${\bf F}(\d)_\bullet$. By the uniqueness of the nonzero maps of the given
degree between the terms of ${\bf F}(\d)_\bullet$, we may identify
$H^0({\cal F}_{\bullet})$ with this complex, proving acyclicity as required for part a).

\medskip

{\bf Case 2)}  $e_m>1$. In this case, Part b) of the Corollary to Bott's
Theorem shows that
each
$H^0({\cal F}_i)$
has an equivariant free resolution of length 1. 
From that Corollary we see moreover that the resolution 
takes the form
$$
0\to {{\bf F}}(\d')_i
 \to {\bf F}(\d)_i \to
  H^0({\cal F}_i) \to 0
$$
for $i=0,\dots,m-1$, 
where $\d'$ is given by
$$
1, d_1+1, \dots, d_{m-1}+1, d_m
$$
corresponding to the sequence
$$
{\bf e}':= (1,e_1 ,\ldots ,e_{m-1}, e_m-1)
$$
and we have simplified the notation by writing
${\bf F}(\d')$ and 
$ {\bf F}(\d)$
instead of 
${\bf F}(\d')(E)$ and 
$ {\bf F}(\d)(E)$.

Because ${\bf F}(\d)(E)_i$ is a free $A$-module generated 
by a representation, and everything splits as $\GL(E)$-modules,
we can lift the differential on $H^0{\cal F}_\bullet$ to get the
following commutative diagram,
where each column is exact. $$
\diagram [small]
  &&0&&0&&&&0\cr
  && \dTo && \dTo &&&& \dTo \cr
0&\rTo&{\bf F}(\d')_m &\rTo&  {\bf F}(\d')_{m-1} &\rTo& \cdots
&\rTo& {\bf F}(\d')_{0} \cr
  &&\dTo&&\dTo&&\cdots&&\dTo                                        \cr
0&\rTo&{\bf F}(\d')_m= {\bf F}(\d)_m&\rTo&  {\bf F}(\d)_{m-1} &\rTo& \cdots
&\rTo& {\bf F}(\d)_{0}   \cr
  &&\dTo&&\dTo&&\cdots&&\dTo                                          \cr
0&\rTo&0                                   &\rTo&     H^0{\cal F}_{m-1}&\rTo& \cdots
&\rTo&H^0 {\cal F}_{0} \cr
  && \dTo && \dTo &&&& \dTo \cr
  &&0&&0&&&&0\cr
\enddiagram
$$
The horizontal maps are constructed
simply to be $\GL(E)$-equivariant
and make the diagram commute, except that we take the upper left
map
${\bf F}(\d')_m \to  {\bf F}(\d')_{m-1}$
to be the map coming from the complex
${\bf F}(\d')_\bullet$,
we take the upper left vertical map
 ${\bf F}(\d')_m\to {\bf F}(\d')_m= {\bf F}(\d)_m$
 to be the equality (so that the left-most column is also exact),
 and we take the horizontal map
 ${\bf F}(\d)_m\to {\bf F}(\d)_{m-1}$
 to be the composite of the two maps above it, assuring the commutativity
 of the upper left-hand square of the diagram. 

We will now prove that the two upper horizontal rows are in fact the complexes
${\bf F}(\d')_\bullet$ and ${\bf F}(\d)_\bullet$. As defined, the
left-hand map in the upper row is the right map, and the map
below it in the middle row is at least nonzero.
From the
uniqueness statement in the definition of the differentials  of our
complexes, all the other horizontal maps will be the correct ones as soon as
we know that they are all nonzero.

Let $\cal G_\bullet$ be the total complex of the double complex
made from the two upper rows of the diagram, so that 
$\cal G_\bullet$ is a resolution of the same module
$H^0(M_{\cal Q})$ as that resolved by
$H^0({\cal F}_{\bullet})$ 
The last vertical
map ${\bf F}(\d')_m\to {\bf F}(\d)_m$ defines
a quotient complex of $\cal G_\bullet$, and is an isomorphism. We may
take the kernel of this quotient map, it has the same homology as $\cal G_\bullet$, 
arriving at a complex
$$
{{\cal G}^\prime_\bullet}:\quad
0\to {{\bf F}}(\d')_{m-1}\to\cdots\to {\bf F}(\d)_{0}
$$
of length $m$ that is, once again, a resolution of the 
module $H^0({M}_{\cal Q})$. 

First, we note that the complex ${{\cal G}^\prime_\bullet}$ is
graded, with degree 0 differentials, if we give the generators of each
 ${\bf F}(\d)_i$ the degree $d_i$ as in the definition of ${\bf F}(\d)_\bullet$,
 and similarly for the  ${\bf F}(\d')_i$. This is because the unique occurence
 of the representation $S_{\alpha(\d,i)}$ that generates 
 ${\bf F}(\d)_i$, in ${\bf F}(\d)_{i-1}$ is in $S_{e_i}\otimes S_{\alpha(\d,i-1)}$, and
 similarly for the ${\bf F}(\d')_i$.
  It follows that all the maps in the resolution $\cal G^\prime_\bullet$
 are given by matrices of elements of positive degree in $A$; that
 is, the resolution $\cal G^\prime_\bullet$ is minimal. From this minimality
 it follows that for each $i$ the constructed map ${\bf F}(\d)_i\to {\bf F}(\d)_{i-1}$ is
  nonzero; for if it vanished then by exactness
 ${\bf F}(\d)_i$ would be in the image of ${\bf F}(\d)_{i+1}\oplus {\bf F}(\d'_i)$,
 which is impossible.
 
 Since the $M_{\cal Q}(\d)$ is the direct sum of finitely many
 Schur functors applied to $\cal Q$, Bott's Theorem tells us that the
 cohomology module 
 $H^0(M_{\cal Q}(\d))$ is a direct sum of finitely many representations,
 each a Schur functor of $E$, and is thus finite-dimensional as
 a vector space. It follows that the dual of ${\cal G}^\prime_\bullet$
 is also a minimal free resolution of an $A$-module of finite
 length. The dual argument to that just given shows that all
 the maps ${\bf F}(\d')_i\to {\bf F}(\d')_{i-1}$ are nonzero as well.
 
 We have now proven the existence of a short exact sequence of
 complexes
 $$
 0\to {\bf F}(\d')_\bullet \to {\bf F}(\d)_\bullet \to H^0({\cal F_\bullet})\to 0.
 $$
We know that the complex
$H^0({\cal F_\bullet})$ is acyclic.
Since $d'_m-d'_0=d_m-d_0-1$, our second induction shows that the 
complex 
${\bf F}(\d')_\bullet$ is acyclic as well. From the long exact sequence
in homology associated to the short exact sequence of complexes, we
see that 
${\bf F}(\d)_\bullet$
is acyclic too, and the proof is done. \hfill $\square$
\medskip

\noindent{\bf Example 3.3} Take $\d= (0,2,5,7,8)$
so that $\e= (0,2,3,2,1)$. Our Young diagram is
$$\Young (0,0,3,4|0,2,3|0,2|1,2|1).$$
The complex ${\bf F}(\d)_\bullet$ has terms
$$(3,1,0,0)\leftarrow (5,1,0,0)\leftarrow (5,4,0,0)\leftarrow (5,4,2,0)\leftarrow (5,4,2,1)$$
where we write $\lambda$ instead of $S_\lambda E\otimes A$.
To get the acyclicity by induction we start with the complex of sheaves on the projective space with the terms
$$(3,1,0)\leftarrow (5,1,0)\leftarrow (5,4,0)\leftarrow (5,4,2)$$
where $\lambda$ is the shorthand for $S_\lambda {\cal Q}\otimes {\cal B}$. Taking modules of sections we get $A$-modules with free resolutions (written as columns)
$$\matrix{(3,1,0,0)&\leftarrow& (5,1,0,0)&\leftarrow &(5,4,0,0)&\leftarrow& (5,4,2,0)\cr
\uparrow&&\uparrow&&\uparrow&&\uparrow \cr
0&\leftarrow&0&\leftarrow&0&\leftarrow&(5,4,2,1)}$$
where again we write $\lambda$ instead of $S_\lambda E\otimes A$. So the mapping cone is the required complex.

\medskip
\noindent{\bf Example 3.4} Take $\d= (0,2,5,6,8)$,
so that $\e= (0,2,3,1,2)$. Our Young diagram is
$$\Young (0,0,0,4|0,2,3,4|0,2|1,2|1).$$
The complex ${\bf F}(\d)_\bullet$ has terms
$$(3,1,1,0)\leftarrow (5,1,1,0)\leftarrow (5,4,1,0)\leftarrow (5,4,2,0)\leftarrow (5,4,2,2)$$
where we write $\lambda$ instead of $S_\lambda E\otimes A$.
To get the acyclicity by induction we start with the complex of sheaves on the projective space with the terms
$$(3,1,1)\leftarrow (5,1,1)\leftarrow (5,4,1)\leftarrow (5,4,2)$$
where $\lambda$ is the shorthand for $S_\lambda {\cal Q}\otimes {\cal B}$. Taking $A$-modules of sections we get $A$-modules with free resolutions (written as columns)
$$\matrix{(3,1,1,0)&\leftarrow& (5,1,1,0)&\leftarrow &(5,4,1,0)&\leftarrow& (5,4,2,0)\cr
\uparrow&&\uparrow&&\uparrow&&\uparrow \cr
(3,1,1,1)&\leftarrow&(5,1,1,1)&\leftarrow&(5,4,1,1)&\leftarrow&(5,4,2,1)}$$
where again we write $\lambda$ instead of $S_\lambda E\otimes A$. 
The first row is the required complex ${\bf F}(\d)_\bullet$ without the last term $(5,4,2,2)$.
But the second row is the complex ${\bf F}(1,3,6,7,8)_\bullet$ without the last term. It corresponds to the Young diagram
$$\Young (0,0,0,0|0,2,3,4|0,2|1,2|1)$$
where the row of zeros is added because $e_0 =1$. Now we notice that the last missing term of this complex is also $(5,4,2,2)$, 
which proves that the homology of the top row is isomorphic to this free $A$-module, and this concludes the proof.

\medskip
\noindent{\bf Remarks} The ranks of the modules in a pure resolution are easy to
calculate from the Herzog-K\"uhl equations; see Section 2.1 of Boij-S\"oderberg [2006]).
In the case of complexes ${\bf F}(\d)_\bullet$ these formulas
are special cases of the Weyl dimension formula. The multiplicative form of these formulas was one of the motivations for looking at Schur functors.

Similarly, it is a standard result that for a graded module of finite length over a polynomial ring, the dimension of the socle equals the Cohen-Macaulay type, 
that is, the rank of the last syzygy module. 
The representation in the highest degree of $M(\d)$ corresponds to the 
partition we get from $\alpha (\d,0)$ by adding one box to each of the first $\alpha(\d,1)_1-1$ rows. 
It is amusing to see that this is the partition $\alpha(\d,m)$ with the first row of length $m$ removed.
So the socle of $M(\d)$ is the representation in the highest degree, as it should since the last
term in the resolution is pure.

We finish this section by analyzing some of the features of the complexes we constructed.

\proclaim Proposition 3.5. Let us assume the sequence $\e$ is symmetric, i.e. $e_i =e_{m+1-i}$ for $i=1,\ldots ,m$. Then the complex $\FF(\d )_\bullet$ is self-dual,
i.e. $\FF(\d )_\bullet^* =\FF(\d )_\bullet$. This duality is $GL(E)$-equivariant.

\noindent{\sl Proof.} Assume that the sequence $\e$ is symmetric. Let $\lambda =\alpha (\d, 0)$. Consider the rectangular partition
$\mu:= (\lambda_1 +e_1)^m$. It is clear that the partitions $\alpha (\d ,i)$ and $\alpha (\d ,m-i)$ are complementary with respect to this rectangle.
This means we have a $GL(E)$-equivariant nondegenerate pairing
$$S_{\alpha (\d ,i)}E\otimes S_{\alpha (\d ,m-i)}E\rightarrow (\bigwedge^m E)^{\otimes(\lambda_1 +e_1)} .$$
This can be extended to the isomorphism of complexes we claim. 
 \hfill $\square$ \medskip

\noindent{\bf Example 3.6} Let us take $m=5$ and $\e = (2,3,1,3,2)$.
The corresponding picture is
$$\Young (0,0,0,0,5|0,0,0,4,5|0,0,0,4,|0,2,3,4,|0,2,,,|1,2,,,|1,,,,).$$
The partitions $\alpha(\d,i)$ for $i = 0,\ldots, 5$ correspond to the
boxes with entries $\leq i$.
Considering the boxes with entries $\geq i+1$ (letting the empty boxes have entry $6$),
if we turn the rectangle $180^\circ$ this is the partition $\alpha(\d,5-i)$.
 
\bigskip
\noindent{\bf \S 4. Modules Supported on Determinantal Varieties.}
\medskip

In this section we describe another way of constructing
a Cohen-Macaulay module whose free resolution is
pure, with given degree shifts $\d$. These
 modules are supported on the degeneracy locus
 of a generic map of free modules $G\to F$, and 
 equivariant for $\GL(F)\times \GL(G)$. 
 (Of course one can derive non-equivariant
artinian modules from them by reducing modulo a general sequence of linear
forms, at least in the case where the ground field is infinite.)
This family of resolutions
 generalizes the ones
 described by Kirby [1974] and Buchsbaum and 
Eisenbud [1975] 
(see Eisenbud [1995] Appendix A2.6 for an exposition) and re-interpreted by 
Weyman ([2003] exercises 37-39, chapter 6), though in the special case treated by those
authors the resolutions work in arbitrary characteristic, while the
method used here to obtain the generalization depends on 
characteristic 0.

With notation as in Section 1,
 we fix the strictly increasing sequence $\d=(d_0,\dots, d_s)$
 and its sequence of differences ${\bf e}=(e_0=d_0 ,e_1=d_1-d_0 ,\ldots )$. Take two vector spaces $F$, $G$, with $\dim(F)= 1+\sum_{i=1}^s(e_i -1)$ and
$\dim(G)=\dim(F)+s-1$. Let $B$ be the polynomial ring
$B= \Sym (F\otimes G^*)$. Consider the Grassmannian 
$\Grass(1,F)$ of lines in $F$ (this is just a projective space),
with tautological sequence
$$
0\rightarrow {\cal R}\rightarrow F\otimes {\cal O}_{\Grass(1, F)}\rightarrow {\cal Q}\rightarrow 0,
$$
so that $\cal Q$ is a bundle of rank $\dim(F)-1$.

Consider the incidence variety 
$$Z=\lbrace (\phi ,{\cal R})\in Hom(F, G)\times Grass(1,F)\ |\ \phi |_{\cal R}=0\rbrace .$$
This is one of the desingularizations of the determinantal variety defined by the maximal minors of the generic matrix $\phi$, 
denoted in the Section 6.5 of Weyman [2003] by $Z^{(1)}_{s-1}$.

Consider the partition 
$\lambda (\d)= ((s-1)^{e_s-1}, (s-2)^{e_{s-1}-1},\ldots, 0^{e_1 -1})$, and let
${\cal N}(\d)$ be the sheaf 
${\cal N}(\d)= S_{\lambda (\d)}{\cal Q}\otimes {\cal O}_Z$. To describe the second family of complexes we set
$$
\gamma(\d, i):= 
((s-1)^{e_s-1}, (s-2)^{e_{s-1}-1},\ldots, i^{e_{i+1}-1}, i^{e_i},(i-1)^{e_{i-1}-1},\ldots ,1^{e_1 -1}).
$$
The partition $\gamma(\d,i)$ is  conjugate to 
the partition $\alpha(\d,i)$ defined in the introduction.

\proclaim Theorem 4.1.
$H^i (Z, {\cal N}(\d ))=0$ for $i>0$, and $H^0 (Z, {\cal N}(\d ))$ has 
a pure $GL(F)\times GL(G)$-equivariant minimal
resolution ${\bf H}(\d )_\bullet$ of type $\d$,
 with terms 
$${\bf H}(\d)_i = 
S_{\gamma(\d,i)}F\otimes\bigwedge^{d_i-d_0}G^*\otimes 
B(-d_i+d_0).
$$

\noindent{\sl Proof.}
Let $p:Z\to Grass(1,F)$ be the projection map. Because $p$ is an affine map it 
suffices for the first statement to show that $H^i p_*({\cal N}(\d))=0$ for $i>0$.
However, 
$$
p_*({\cal N}(\d))=S_{\lambda(\d)}{\cal Q}\otimes Sym({\cal Q}\otimes G^*).
$$
Since this does not involve $\cal R$, it has no higher cohomology.
 
To prove the second statement,
we apply the Basic Theorem (5.1.2) from Weyman [2003] to the sheaf 
${\cal N}(\d)$. In the notation of that result, we set $\xi = {\cal R}\otimes G^*$, $\eta = {\cal Q}\otimes G^*$ and ${\cal V}=S_{\lambda (\d)}{\cal Q}$.
We get a complex ${\bf F}(S_{\lambda (\d)}{\cal Q})_\bullet$ which is our ${\bf H}(\d )_\bullet$, which is a resolution of 
$H^0 p_*({\cal N}(\d))$
because the higher cohomology
$H^i p_*({\cal N}(\d))$ vanishes.

The direct calculation of the cohomology groups using Bott's theorem (Weyman [2003],
 (4.1.9)), dualized using exercise 18 b), p.83),
gives the terms of our complex. More precisely, the calculation comes down to applying Bott's Theorem to the weights
$$((s-1)^{e_s-1}, (s-2)^{e_{s-1}-1},\ldots, 0^{e_1 -1},u)$$
for $0\le u\le dim(G)$. The partition $\gamma(\d,i)$ comes from the term with
$u=d_i-d_1+1+i$.
\hfill$\square$\medskip

\noindent{\bf Remark} The above theorem is the specialization of Theorem 0.2 we
get by setting $V = V_1 = F$ and $U = U_1 = G^*$.

\bigskip
\noindent{\bf \S 5. Proofs of Theorems 0.1 and 0.2.}
\medskip
\noindent {\sl Proof of Theorem 0.1.} The differential in $\FF(\d )_\bullet$ is $\gl (V)$-equivariant by definition. The $\ZZ/2$-graded Pieri formula implies $\FF(\d )_\bullet$ is a complex.
By Berele-Regev theory, the $\gl (V)$ action on each homogeneous component of each module
$\FF(\d )_i$ is semi-simple. Thus the homogeneous components of the homology are direct sums of Schur modules $\SS_\mu V$.

Now assume that for some $\lambda$ and some $(m,n)$ the complex $\FF(\d )_\bullet$ is not acyclic, so some module $\SS_\mu V$ consists of cycles that are not boundaries. But then, for a bigger dimension vector $(m, n)$, such that $S_\mu (V_0 )\ne 0$, the same Schur module $\SS_\mu V$ is in the homology.
If $m$ is large enough,
this Schur module contains elements in some weight that does not involve $V_1$, that is,
elements defined by tableaux that contain only basis elements from $V_0$.

Since the differential of $\FF (\d )_\bullet$ preserves
the weight space decomposition, this implies that the complex
$\FF (\d )_\bullet$ for $V_0$ is not acyclic. This contradicts Theorem 3.2,
proving the Theorem.
 \hfill $\square$ \medskip

\bigskip

\noindent{\sl Proof of Theorem 0.2.} The proof of Theorem 0.2 follows the same
outline as the proof of Theorem 0.1, except that we use subspaces
$U_1\subset U$ and $V_1\subset V$, and Theorem 4.1, where we set
$F=V_1$ and $G^*=U_1$. \hfill $\square$ \medskip

\medskip
\noindent{\bf Remark} Consider Theorem 0.1 in the odd case, that is where
 $m=0$, and set $F=V_1$. 
Then $R=\bigwedge^\bullet (F)$ and the complex  $\FF(\d)_\bullet$ has terms
$$
\FF(\d )_i = S_{\alpha(\d, i)'}E\otimes R(-e_1-\ldots -e_i).
$$
In this way we get truncations of the Tate resolutions constructed by Fl\o ystad [2004].  An alternate proof of Theorem 0.1. could be obtained by reducing to the odd case.
The proof in the current paper avoids Bernstein-Gel'fand-Gel'fand duality. 
One can also relate the even and odd cases to each other by means of Schur duality.

\bigskip
\noindent{\bf \S 6. Some open problems and conjectures.}
\medskip

In 6.1--6.5 below we work over the polynomial ring $A$.

\proclaim Conjecture 6.1. Every sufficiently large integral point in the
ray defining the possible pure Betti tables of graded modules of
finite length over a polynomial ring, with a given degree
sequence, is actually the Betti table of the free resolution of 
a Cohen-Macaulay module.

The following examples are expressed in terms of the smallest integral
multiple of the Betti number on a ray of pure resolutions, corresponding
to the given degree sequence $\d$,
which we call the primitive vector of Betti numbers. In all these examples
$m=3$.

\medskip
\noindent
{\bf Example 6.2.} $\d = (0,3,4,7)$.
The primitive vector of Betti numbers is $\beta =(1,7,7,1)$. The
construction of \S 3 gives
Betti numbers $6\beta$, the Eisenbud-Schreyer construction gives
$15\beta$ and the construction from \S 4 gives $50\beta$. Thus all three
are needed in order to conclude the conjecture for this extremal ray.
But in this case we know that the primitive vector is achieved
by the minimal free resolution of the ideal of $6\times 6$
 Pfaffians of a $7\times 7$ skew-symmetric matrix of linear forms; see 
 Buchsbaum and Eisenbud [1975].

\medskip
\noindent
{\bf Example 6.3.} $\d = (0,4,9,13)$.
The primitive vector of Betti numbers is $\beta =(5,13,13,5)$. The
construction of \S 3 gives
Betti numbers $18\beta$, the Eisenbud-Schreyer construction gives
$380\beta$ and the construction of \S 4 gives $9075\beta$. 
These three examples together imply the truth of the 
conjecture on this ray.

\medskip
\noindent
{\bf Example 6.4.} $\d = (0,1,4,6)$. The primitive
vector of Betti numbers is $\beta = (5,8,5,2)$. All three
constructions give $5\beta$. This is the smallest sequence for $n=3$
where we cannot conclude the conjecture using our three constructions.
\medskip
\noindent 
The material presented in this paper raises some other interesting questions.

\smallskip
\noindent
{\bf Problem 6.5.} Equivariant Boij-S\"oderberg conjectures. We use the
notation of Section 3.
Let $\FF_\bullet$ be a $GL(E)$-equivariant acyclic complex of
$A$-modules. Is it in the cone generated by the Betti tables of pure
resolutions constructed in Section 3 in the sense that there
exist $GL(E)$-representations $W$ and $W_1 ,\ldots ,W_s$ and degree
shifts $\d (i)$ (for $1\le i\le s$) such that for each $j=0,\ldots ,m$
we have isomorphisms of $GL(E)$-modules
$$W\otimes \FF_j =\bigoplus_{i=1}^s W_i\otimes \FF(\d (i))_j ?$$
In particular, assuming $\FF_\bullet$ is pure with degree shifts $\d$,
does it mean that $W\otimes \FF_\bullet =W'\otimes \FF(\d )_\bullet$
for some $GL(E)$-modules $W, W'$?

\medskip
\noindent
{\bf Problem 6.6.} $\ZZ/2$-graded Betti tables. Let 
$R= Sym(V_0 )\otimes \bigwedge^\bullet (V_1).$
 Is
the Betti table of every acyclic complex of free $R$-modules in the
cone generated by the Betti tables of pure acyclic complexes ? What
are the facets of the cone generated by Betti tables of pure acyclic
complexes of free modules, is this cone self-dual in some sense ?

\bigskip
\goodbreak

\centerline{\bf References.}

\item{}A. Berele and A. Regev: Hook Young diagrams with applications to combinatorics and to representations of Lie superalgebras. Adv. in Math. 64 (1987) 118-175. MR 88i:20006.
\item{} M.~Boij, J. S\"oderberg, Graded Betti numbers of Cohen-Macaulay modules and the multiplicity conjecture, arXiv:math/0611081.
\smallskip\item{} D.~A.~Buchsbaum, D.~Eisenbud, Generic free resolutions and a family of generically perfect ideals. Advances in Math. 18 (1975), no. 3, 245--301.
\smallskip\item{} D.~A.~Buchsbaum, D.~Eisenbud, Algebra structures for finite free resolutions, and some structure theorems for ideals of codimension $3$.  Amer. J. Math.  99  (1977), no. 3, 447--485.
\smallskip\item{} M. Demazure, A very simple proof of Bott's theorem, Inventiones Mathematicae 34 (1976), 271-272.
\smallskip\item{} D.~Eisenbud, {\it Commutative Algebra with a View Toward Algebraic Geometry}, Graduate Texts in Mathematics, Springer, 1995.
\smallskip\item{} D.~Eisenbud, F.~Schreyer, Betti Numbers of Graded Modules and Cohomology of Vector Bundles, arXiv0712.1843.
\smallskip\item{}D.~Eisenbud, J.~Weyman, Fitting's Lemma for ${\bf Z}/2$-graded modules,
Trans. Am. Math. Soc. 355,  4451-4473 (2003)
\smallskip\item{} G.~Fl\o ystad, 
Exterior algebra resolutions arising from homogeneous bundles. 
Math. Scand. 94 (2004), no. 2, 191--201. 
\smallskip\item{} W.~Fulton, J.~Harris, {\it Representation Theory; the first course}, 
Graduate Texts in Mathematics 129, Springer-Verlag, 1991.
\smallskip\item{} J.~Herzog and M.~K\"uhl, M.
On the Betti numbers of finite pure and linear resolutions.
Comm. Algebra 12 (1984), no. 13-14, 1627--1646. 
\smallskip\item{} D.~Kirby,A sequence of complexes associated with a matrix.
J. London Math. Soc. (2) 7 (1974), 523--530. 
\smallskip\item{}  I.~G.~Macdonald: {\it Symmetric functions and Hall polynomials.}
Second edition. With contributions by A.
Zelevinsky. Oxford Mathematical Monographs,
Oxford University Press, New York, 1995.
\smallskip\item{} C. Peskine, and L.~Szpiro,
Dimension projective finie et cohomologie locale. Applications \`a la d\'emonstration de 
conjectures de M. Auslander, H. Bass et A. Grothendieck.
Inst. des Hautes \'Etudes Sci. Publ. Math. No. 42 (1973), 47--119.
\smallskip\item{}  S.~Sam and J.~Weyman, Pieri Resolutions for Classical Groups. arXiv:0907.4505 .

\smallskip\item{} J.~Weyman, {\it Cohomology of vector bundles and syzygies}, Cambridge University Press,  Cambridge (2003).

\bye